\documentclass[10pt]{article}
\font\smallit=cmti10

\usepackage{amssymb,latexsym,amsmath,epsfig,amsthm} %% Add other packages as necessary

\makeatletter

\renewcommand\section{\@startsection {section}{1}{\z@}
{-30pt \@plus -1ex \@minus -.2ex}
{2.3ex \@plus.2ex}
{\normalfont\normalsize\bfseries\boldmath}}

\renewcommand\subsection{\@startsection{subsection}{2}{\z@}
{-3.25ex\@plus -1ex \@minus -.2ex}
{1.5ex \@plus .2ex}
{\normalfont\normalsize\bfseries\boldmath}}

\renewcommand{\@seccntformat}[1]{\csname the#1\endcsname. }

\makeatother

\newtheorem{theorem}{Theorem}

\theoremstyle{definition}

\begin{document}

\begin{center}
\uppercase{\bf Lower bounds for ranks using Pell equations.}
\vskip 20pt
{\bf P.G. Walsh}\\
{\smallit Department of Mathematics, University of Ottawa, Ottawa, Ontario, Canada}\\
{\tt gwalsh@uottawa.ca}\\
\vskip 10pt
\end{center}

\centerline{\bf Abstract}
\noindent
In this short note, we examine the ranks of a subfamily of curves from a previous paper derived from the existence of solutions to certain Pell equations. We exhibit an abundance of curves of moderately large rank, and using certain well known results from Diophantine analysis, we prove under mild conditions that these curves have rank at least three provided that the constant coefficient of the cubic polynomial defining the curve is sufficiently large.\\


\pagestyle{myheadings}
%\markright{\smalltt INTEGERS: 21 (2021)\hfill}
\thispagestyle{empty}
\baselineskip=12.875pt
\vskip 30pt

\section{Introduction}
In a previous paper \cite{W2}, it was shown that an elliptic curve taking the form
$$y^2=x(x+a)(x+b)+m^2,$$
with $a,b$ distinct non-zero integers, and any sufficiently large integer $m$, has rank at least $2$. This result was motivated by the main result in the seminal paper \cite{BM} of Brown and Myers, which considered the particular case $a=1,b=-1$. Related to their paper is a recent paper by Hatley and Stack \cite{HS}, in which they study the subfamily of curves in \cite{BM} with $m^2$ replaced by $m^6$. This latter modification can be a profitable one from the point of view of increasing the rank of the curve. In particular, if $a$ and $b$ defining the curve above give rise to a certain Pell equation which is solvable, one can show the following, which constitutes the main result of this paper.

\begin{theorem}
Let $a$ and $b$ be distinct coprime integers for which the Pell equation
$$X^2-(a+b)Y^2=-ab$$
is solvable, and let $(X,Y)=(n,m)$ be an integer solution. Assume further that the polynomials $F_{a,b}(x,m)$, $G_{a,b}(x,m)$ and $H_{a,b}(x,m)$ in the proof below are irreducible in $\mathbb{Q}[x,m]$. Then for $m$ sufficiently large, the rank of the curve
$$E: \; \; y^2=x(x+a)(x+b)+m^6 \eqno (1.2)$$
is at least $3$.
\end{theorem}

The polynomials $F_{a,b}(x,m)$, $G_{a,b}(x,m)$ and $H_{a,b}(x,m)$ are irreducible as polynomials in $\mathbb{Q}[x,m,a,b]$, and our ongoing computations indicate that they are irreducible in $\mathbb{Q}[x,m]$ for all distinct positive integers $a$ and $b$. The polynomial $F_{a,b}(x,m)$ is given explicitly in the proof of Theorem 1.1, however $G_{a,b}(x,m)$ and $H_{a,b}(x,m)$ are sufficiently long that we refer the reader to \cite{W3} to access them.\\

At this point we exhibit a few relatively small examples in order to give the reader a sense of how useful this construction can be for finding curves with rank of moderate size, although we make absolutely no claims in the direction of breaking any records.\\

\noindent {\bf Example 1.1} Let $(a,b)=(1,2)$. The Pell equation $X^2-3Y^2=-2$
has the solution $(X,Y)=(1,1)$, which we write as $\alpha = 1+\sqrt{3}$, giving $m=1$, and in this case, the curve defined in Theorem 1.1 has rank $1$. However, as we multiply $\alpha$ by powers of the fundamental unit $\epsilon_3=2+\sqrt{3}$ to get more solutions to the above Pell equation, things change in our favour in quite a hurry.\\

\noindent $k=0$, $\alpha \epsilon_3^k=1+\sqrt{3}$, $E: y^2=x^3+3x^2+2x+1$, $r=1$.\\
$k=1$, $\alpha \epsilon_3^k=5+3\sqrt{3}$, $E: y^2=x^3+3x^2+2x+3^6$, $r=4$.\\
$k=2$, $\alpha \epsilon_3^k=19+11\sqrt{3}$, $E: y^2=x^3+3x^2+2x+11^6$, $r=5$.\\
$k=3$, $\alpha \epsilon_3^k=71+41\sqrt{3}$, $E: y^2=x^3+3x^2+2x+41^6$, $r=7$.\\
$k=4$, $\alpha \epsilon_3^k=265+153\sqrt{3}$, $E: y^2=x^3+3x^2+2x+153^6$, $r=5$.\\
$k=5$, $\alpha \epsilon_3^k=989+571\sqrt{3}$, $E: y^2=x^3+3x^2+2x+571^6$, $r=6$.\\
$k=6$, $\alpha \epsilon_3^k=3691+2131\sqrt{3}$, $E: y^2=x^3+3x^2+2x+2131^6$, $r=6$.\\
$\; $\newline

\noindent {\bf Example 1.2} The Pell equation in the statement of Theorem 1.1 is always solvable when $a=1$, and to exploit the existence of a small fundamental unit, as in the previous example, one can put $b=t^2-2$ so that the discriminant $a+b=t^2-1$. In this way, the first few corresponding solutions to the Pell equations, and their corresponding parametric family of elliptic curves, are as follows.\\

\noindent $k=0$, $\alpha \epsilon^k = 1+\sqrt{t^2-1}$, \newline
$E_0(t): \; y^2=x^3+(t^2-1)x^2+(t^2-2)x+1.$\\

\noindent $k=1$, $\alpha \epsilon^k = t^2+t-1+(t+1)\sqrt{t^2-1}$, \newline
$E_1(t): \; y^2=x^3+(t^2-1)x^2+(t^2-2)x+(t+1)^6.$\\

\noindent $k=2$, $\alpha \epsilon^k = 2t^3+2t^2-2t-1+(2t^2+2t-1)\sqrt{t^2-1}$, \newline
$E_2(t): \; y^2=x^3+(t^2-1)x^2+(t^2-2)x+(2t^2+2t-1)^6.$\\

The reader may wish to attempt to compute the ranks of these curves, although their heights grow extremely rapidly, perhaps requiring developing methods which make use of Mestre-Nagao sums (the reader is referred to \cite{KV} and its references for more on this topic). At the time of writing, the record is $rk(E_1(6001))=9$ sent to the author recently by Andrej Dujella.

\section{Proof of Theorem 1.1}

We now turn our attention to the proof of Theorem 1.1. It is already known from \cite{W2} that $E$ is $2$-torsion free, and that $P=(-a,m^3)$ and $Q=(-b,m^3)$ are independent points on E. Let $R$ denote the point $R=(-m^2,mn)$, which is on $E$ because of the fact that $(n,m)$ is an integer solution to the Pell equation in the statement of the theorem. In order to show that $R$ is a third independent point on E, we apply Proposition 1.5 of \cite{HS} by showing that none of $R$, $P+R$, $Q+R$ and $P+Q+R$ are doubles of a point on $E$. The polynomials $F$, $G$, and $H$ arise from each of these cases respectively (the case $P+R$ is the same as $Q+R$).\\

Let $(x,y)$ denote a point on $E$, and assume that $2(x,y)=R=(-m^2,mn)$. Using the doubling formulae on $E$ from p.54 of \cite{Sil}, equating the corresponding quantities in the $x$ coordinates, and translating $x$ by $m^2$, one obtains the equation $F_{a,b}(x,m)=0$, where $F_{a,b}(x,m)$ is given by
$$x^4 + (4am^2 + 4bm^2 -6m^4 - 2ab)x^2 + (8abm^2 - 8am^4 - 8bm^4)x + (9m^8 - 6m^4ab + a^2b^2).$$
In order to derive the desired result, one requires that certain properties hold with regard to $F_{a,b}(x,m)$. Firstly, the curve defined by $F_{a,b}(x,m)=0$ is singular, but can be desingularized by the map $x \rightarrow xm$, and the resulting curve has genus equal to 5. Therefore, if it were known that $F_{a,b}(x,m)$ is irreducible as a polynomial in $\mathbb{Q}[x,m]$ for every pair $a,b$ being considered in the statement of the theorem, the proof that $R$ is not in $2E$ would be complete by Faltings theorem.\\

In order to show that none of $P+R,Q+R$ and $P+Q+R$ are in $2E$, we follow the same approach. A computation almost identical to that above, equating the $x$-coordinate of $2(x,y)$ and that of $P+R$, yields a polynomial $G_{a,b}(x,m)$ of degree $8$ in $x$. Thus, a solution to $2(x,y)=P+R$ gives rise to an integer solution to $G_{a,b}(x,m)=0$. Similar to $F_{a,b}(x,m)$, the polynomial $G_{a,b}(x,m)$ is irreducible in $\mathbb{Q}[x,m,a,b]$ when regarded as a polynomial in four variables, and has been computationally verified to be irreducible in $\mathbb{Q}[x,m]$ for all distinct positive integers $1 < a < b$ up to $10^3$. For $a=1$, $G$ is divisible by $m^2-1$, and so in this case, $G$ is replaced by the cofactor. The weighted sum of highest order terms of $G_{a,b}(x,m)$ satisfies the reducibility hypothesis of Runge's theorem on Diophantine equations (see \cite{W1}). Therefore, under the irreducibility condition, one can assert that $m$ is bounded effectively in terms of $a$ and $b$. The interested reader can access $G_{a,b}(x,m)$ directly from \cite{W3}.\\

Finally, for the case $P+Q+R$, one can follow the identical procedure using the equation $2(x,y)=P+Q+R$. This results in an integer solution to an equation of the form $H_{a,b}(x,m)=0$ in which the polynomial $H_{a,b}(x,m)$ is irreducible in $\mathbb{Q}[x,m,a,b]$, and has been verified to be irreducible in $\mathbb{Q}[x,m]$ for all pairs of distinct positive integers $(a,b)$ up to $10^3$. The curve defined by $H_{a,b}(x,m)=0$ is of positive genus, and hence the argument given above for $F_{a,b}(x.m)$ shows that the equation $H_{a,b}(x,m)=0$ is not solvable in integers $(x,m)$ for $m$ large, which completes the proof of the theorem.\\

\noindent {\bf Remark.} Each of the parametric curves $E_i(t)$ ($-10 \le i \le 10$) have been verified using Magma to have rank at least 3 over $Q(t)$ by showing that $P$, $Q$ and $R$ are linearly independent over the function field. Therefore, Silverman's Specialization Theorem (see Theorem 20.3 in \cite{Sil}) provides an effective proof that the curves in these families have rank at least $3$ for $t$ sufficiently large.\\

\noindent {\bf Acknowledgements} The author gratefully acknowledges the communications with Andrej Dujella, Noam Elkies and Joe Silverman.

\end{document}